\documentclass{article}[12pt]
\usepackage{amsfonts,amsmath}

\usepackage{color}

\newtheorem{theorem}{Theorem}
\newtheorem{proposition}{Proposition}
\newtheorem{lemma}{Lemma}

\def\R{{\mathbb R}}  
\def\C{{\mathbb C}}

\begin{document}

\title{Darboux--Moutard transformations and Poincare--Steklov operators
\thanks{The work was supported by the French--Russian grant (RFBR 17-51-150001 NCNI$\_$a/PRC 1545 CNRS/RFBR) and done during the visit of the second author (I.A.T.) to Centre de Math\'ematiques Appliqu\'ees of
\'Ecole Polytechnique.}}

\author{R.G. Novikov \thanks{CNRS (UMR 7641), Centre de Math\'ematiques
Appliqu\'ees, \'Ecole Polytechnique, 91128 Palaiseau, France; e-mail:
novikov@cmap.polytechnique.fr} \and I.A. Taimanov
\thanks{Sobolev Institute of Mathematics, 630090 Novosibirsk,
Russia, and Novosibirsk State University, 630090 Novosibirsk, Russia; e-mail: taimanov@math.nsc.ru}}

\date{}

\maketitle

\hfill{To S.P. Novikov on the 80th birthday}

\begin{abstract}
Formulas relating Poincare--Steklov operators for Schr\"odinger equations related by 
Darboux--Moutard transformations are derived. They can be used for testing algorithms of reconstruction of
the potential from measurements at the boundary.
\end{abstract}

Investigations of inverse problems for two-dimensional Schr\"odinger operators at a given energy level were initiated within the framework of the theory of solitons by S.P. Novikov and his scientific school. In particular, in \cite{1} the spectral data were introduced for the two-dimensional periodic Schr\"odinger operator, in a magnetic field, which is finite-gap at one energy level, and the inverse problem of reconstructing the operator from these data was solved; in \cite{2,3}, in terms of these data, potential operators were singled out and there were derived the evolutionary equations (the Novikov--Veselov equations) 
which preserve this class of operators and their spectra at the given energy level, and theta-functional formulas for solving the equations were obtained; in \cite{4} the first results on the inverse scattering problem at the negative energy level 
were obtained and for the first time in the theory of inverse problems the methods of the theory of generalized analytic functions were used.

In this paper we develop an approach to direct and inverse problems, for the two-dimensional  Schr\"odinger operator,
based on the Moutard transformation. Since the Darboux transformation is a one-dimensional reduction of the Moutard transformation, our results on two-dmensional operators from \S 2 are extended to the case of one-dimensional
 Schr\"odinger operators (see \S 3).
 
\section{Preliminary facts}

\subsection{Darboux--Moutard transformations}

The Moutard transformation constructs from solutions of a second order equation of the type
$$
\partial_v \partial_w \varphi + u \varphi =0
$$
solutions $\vartheta$ of another second order equation of the same type
$$
\partial_v \partial_w \theta  +\widetilde{u} \theta =0
$$
and is determined by a solution of the initial equation \cite{Moutard}.
We consider its special reduction for which the variables $v$ and $w$ are complex--conjugate
and after renormalizations the second order equation reduces to the  Schr\"odinger equation
\begin{equation}
\label{1}
H \varphi = -\Delta \varphi + u \varphi = 0,
\end{equation}
where
$$
\Delta = \frac{\partial^2}{\partial x^2} + \frac{\partial^2}{\partial y^2}
$$
is the Laplace operator on the two-plane $\R^2$.
Let us take a solution $\omega$ of (\ref{1})
and construct by using this solution the new Schr\"odinger operator
$$
\widetilde{H} = -\Delta + \widetilde{u}
$$
with the potential
\begin{equation}
\label{2}
\widetilde{u} =  u -2\Delta \log \omega =  - u +
2\frac{\omega_x^2+\omega_y^2}{\omega^2}.
\end{equation}

Is is easy to show by straightforward computations that if a function
$\varphi$ satisfies (\ref{1}), then the function $\theta$, which is determined by the relations
\begin{equation}
\label{3}
(\omega \theta)_x = -\omega^2
\left(\frac{\varphi}{\omega}\right)_y, \ \ \ (\omega \theta)_y =
\omega^2 \left(\frac{\varphi}{\omega}\right)_x,
\end{equation}
satisfies the equation
\begin{equation}
\label{4}
\widetilde{H}\theta = -\Delta \theta + \widetilde{u}\theta = 0,
\end{equation}
obtained from (\ref{1}) by the Moutard transformation determined by the initial solution
$\omega$.
We note that the relations (\ref{3}) determine $\theta$ up to summands of the form
$\frac{C}{\omega}, C = \mathrm{const}$,
the function $\omega^{-1} =\frac{1}{\omega}$ satisfies the equation
$$
\widetilde{H}\frac{1}{\omega} = 0
$$
and determines the inverse Moutard transformation from $\widetilde{H}$ to $H$:
$$
u \ \ \stackrel{\omega}{\longrightarrow} \ \ \widetilde{u} = u  - 2\Delta \log \omega \ \ \stackrel{1/\omega}{\longrightarrow} \ \
u = \widetilde{u} - 2\Delta \log \frac{1}{\omega} = \widetilde{u} + 2\Delta \log \omega.
$$

For a potential
$$
u = u(x)
$$
which depends on one variable the Moutard transformation reduces to the Darboux transformation
\cite{Darboux}.
Let $H_1$ be a one-dimensional Schr\"odinger operator
$$
H_1 = -\frac{d^2}{dx^2} + u
$$
and $\omega_1$ be its eigenfunction:
$$
H_1 \omega_1 = E \omega_1, \ \ \ E_0 = \kappa^2.
$$
The Moutard transformation, of the two-dimensional Schr\"odinger operator, deter\-mi\-ned by the solution
$$
\omega (x,y) = e^{\kappa y} \omega_1(x)
$$
of (\ref{1}), maps the potential of the operator into
$$
\widetilde{u}(x)  = u(x)  - 2\Delta \log (e^{\kappa y} \omega_1(x)) = u - 2 \frac{d^2}{dx^2} \log \omega_1(x).
$$
Therewith the corresponding one-dimensional Schr\"odinger operator $H_1$ is trans\-formed into the one-dimensional
operator
$$
\widetilde{H}_1 = -\frac{d^2}{dx^2} + \widetilde{u}.
$$
This transformation is called the Darboux transformation. Its action on eigen\-func\-tions takes the following form.
Let
$$
H_1 \varphi_1 = E \varphi_1, \ \ \ \varphi_1=\varphi_1(x), \ \ E=\mu^2.
$$
We put
$$
\varphi = e^{\mu y}\varphi_1(x).
$$
It is clear that $H \varphi=0$ and, by (\ref{3}),
we get the Moutard transformation of $\varphi$ in the form
\begin{equation}
\label{5}
\theta = e^{\mu y}\theta_1(x),
\end{equation}
were the second of relations (\ref{3}) gives the Darboux transformation of eigen\-func\-tions:
\begin{equation}
\label{6}
\theta_1 = \frac{1}{\mu +\kappa}\left(\frac{d}{dx} - \frac{d \log \omega_1}{dx}\right)\varphi_1,
\end{equation}
$$
\widetilde{H}_1 \theta_1 = -\frac{d^2 \theta_1}{dx^2} + \widetilde{u} \theta_1 = E \theta_1, \ \ \ E=\mu^2,
$$
and the first one gives its inversion
\begin{equation}
\label{7}
\varphi_1 = \frac{1}{\kappa-\mu}\left(\frac{d}{dx} + \frac{d \log \omega_1}{dx}\right)\theta_1.
\end{equation}
In this case, since the image of the transformation is sought in the form (\ref{5}), the transformation of eigenfunctions becomes single-valued, because the addition of summands of multiple $ \frac{1}{\omega} $ does not preserve the form (\ref{5})
(for $\mu \neq \kappa $).

The Darboux transformation was repeatedly used and was often rediscovered
(see, for example, \cite{Crum}) for solving problems of mathematical physics and the spectral theory \cite{MS, M} (see also the review \cite{TT2010}).

The Moutard transformation and its extension
for solutions of the Novikov--Veselov equations \cite{HLL,TT2007} in recent years has been applied to the construction of the first examples of two-dimensional Schr\"odinger operators with fast decaying potential and with a nontrivial kernel \cite{TT2007} and
blowing-up solutions of the Novikov--Veselov equation with regular initial data \cite{TT2008a, TT2008} (see the numerical analysis of the negative discrete spectrum and its dynamics for these examples in \cite{AT}), to the construction of explicit examples of two-dimensional potentials of Wigner--von Neumann type \cite{NTT2014}. In \cite{TT2013}, the action of the Moutard transformation on the Faddeev eigenfunctions at the zero energy level was described, and in\cite{NT2013}, using this transformation, Faddeev's eigenfunctions at the zero energy level for multipoint delta-like potentials were found.

In \cite{T2015a, MT2017}, it was established a relation of the generalized Moutard transformation
for two-dimensional Dirac operators \cite{C} to the con\-for\-mal geometry of surfaces in three- and four-dimensional spaces, and with this were constructed blowing-up solutions of the modified Novikov--Veselov equation
with regular initial data \cite{T2015b,T2015c}.
A generalization of the Moutard transformation to the case of generalized analytic functions, in particular, gave an approach to constructing the theory of generalized analytic functions with contour poles \cite{NG1, NG2, NG3, NT2016} and also allowed to construct a Moutard--type trans\-for\-mation for the conductivity equation \cite{NG4}.

\subsection{Poincare--Steklov operators}

Let in the domain $D$ with the boundary $\partial D$ there is given an elliptic differential equation
\begin{equation}
\label{8}
L \psi = E\psi.
\end{equation}
We single out two boundary conditions, each of which, as a rule, completely determines the solution of the equation.
Then the Poincare-Steklov operator, by definition, takes the value of one
boundary condition into the value of another condition.

Let us consider the most well-known particular cases of Poincare-Steklov operators.
Let the elliptic equation (\ref{8}) be given by a linear differential
second-order expression $L$. Then

\begin{enumerate}
\item
if $ E $ is not an eigenvalue of the problem (\ref{8})
with the Dirichlet condition $\psi|_{\partial D} = 0$, then
the boundary data $\psi|_{\partial D}$
determine the solution of (\ref{8}) uniquely and we define the DN
(Dirichlet-to-Neumann) operator, which takes the values of $\psi$ on the boundary to the values of the derivatives
of  $\psi$ along the exterior normal $\nu$ to the boundary (the data of the Neumann problem):
\begin{equation}
\label{9}
DN: \ \psi|_{\partial D} \longrightarrow \frac{\partial \psi}{\partial \nu}|_{\partial D};
\end{equation}

\item
if $E$ is not an eigenvalue of the problem (\ref{8})
with the Neumann condition $\frac{\partial\psi}{\partial\nu}|_{\partial D} = 0$, then
the ND (Neumann-to-Dirichlet) operator is defined:
\begin{equation}
\label{10}
ND: \ \frac{\partial \psi}{\partial \nu}|_{\partial D} \longrightarrow \psi|_{\partial D};
\end{equation}

\item
the above operators are special cases of the RR (Robin-to-Robin) operator, which in the general case relates
mixed boundary conditions (the Robin conditions). If $E$ is not an eigenvalue of the problem (\ref{8})
with the boundary condition
$$
\left(\cos \alpha \, \psi - \sin \alpha \,\frac{\partial \psi}{\partial \nu}\right)|_{\partial D} =0,
$$
then the RR operator maps the boundary data
$$
\Gamma_\alpha \psi  = \left(\cos \alpha \,\psi - \sin \alpha \,\frac{\partial \psi}{\partial \nu}\right)|_{\partial D}
$$
to the boundary data $\Gamma_{\alpha - \pi/2}\psi$:
$$
RR: \ \Gamma_\alpha \psi \longrightarrow \Gamma_{\alpha-\pi/2} \psi.
$$
As particular cases, we get the DN operator for
$\alpha =0$ and the ND operator for $\alpha = \frac{\pi}{2}$.
\end{enumerate}

\section{The action of the Moutard transformation on Poincare-Steklov operators}

We assume that equation (\ref{1}) holds in a bounded simply-connected two-dimen\-sio\-nal domain 
$D \subset \R^2$ with a smooth boundary $\partial D$ and that $u$ is a regular function on $D \cup \partial D$.

Let $\widetilde{H}$ and $\widetilde{u}$ be the operator and the potential from (\ref{4}),
i.e., the Moutard transformations of $H$ and $u$, determined by a solution $\omega$
of (\ref{1}) via formulas (\ref{2}).

For equation (\ref{1}) we consider the operator $$\Phi_u=DN$$ of the form (\ref{9}) with $E=0$ and the operator $$\Phi^{-1}_u = ND$$ of the form (\ref{10}) with $E=0$, where $u$ is the potential from (\ref{1}).
Let $\sigma_D(H)$ and $\sigma_N(H)$ denote the spectra of the operators defined by $H = -\Delta+u$ and
the Dirichlet and Neumann boundary conditions, respectively.

For solutions $\varphi$ of (\ref{1}) we consider also the following boundary data on $\partial D$:
\begin{equation}
\label{11}
\Gamma^\tau_\omega  \varphi= \omega \left(\frac{\varphi}{\omega}\right)_\tau = \varphi_\tau -
\frac{\omega_\tau}{\omega} \varphi
\end{equation}
and
\begin{equation}
\label{12}
\Gamma^\nu_\omega \varphi = \omega \left(\frac{\varphi}{\omega}\right)_\nu = \varphi_\nu -
\frac{\omega_\nu}{\omega} \varphi,
\end{equation}
where $\omega$ is the fixed solution of (\ref{1}) in $D$, $\nu$ is the outer normal to $\partial D$,
$\tau$ is the path-length parameter on $\partial D$, which grows in the direction of $\nu^\perp = (-\nu_2, \nu_1)$ with $\nu = (\nu_1,\nu_2)$; the lower indices $\tau$ and $\nu$ denote the derivations in $\tau$ 
and along the normal $\nu$.

For simplicity, we assume that $\omega$ has no zeroes on $\partial D$.

The important observation is that relations (\ref{3}) on $\partial D$ can be rewritten as
 \begin{equation}
 \label{13}
 \Gamma^\nu_{\omega^{-1}} \theta = - \Gamma^\tau_\omega \varphi,
 \end{equation}
 \begin{equation}
 \label{14}
 \Gamma^\tau_{\omega^{-1}}\theta = \Gamma^\nu_\omega \varphi.
 \end{equation}

\begin{lemma}
Assuming that $0 \notin \sigma_D(H)$, the following formulas hold:
$$
\mathrm{Ker}\, \Gamma^\tau_\omega = \{c \omega \ : \ c\in \C\},
$$
$$
\mathrm{Im}\, \Gamma^\tau_\omega = \{ f \ : \ \int_{\partial D} \omega^{-1} f \,d\tau =0\},
$$
where $\Gamma^\tau_\omega$ is considered as an operator on solutions of (\ref{1}).
\end{lemma}

Lemma 1 follows from the definition of $\Gamma^\tau_\varphi$ by (\ref{11}) and the fact that under our assumptions the Dirichlet problem for (\ref{1}) is uniquely solvable.

\begin{lemma}
Assuming that $0 \notin \sigma_D(\widetilde{H})$, we have
\begin{equation}
\label{15}
\mathrm{Ker}\, \Gamma^\nu_\omega = \{c \omega \ : \ c \in \C\},
\end{equation}
$$
\mathrm{Im}\, \Gamma^\nu_\omega = \{ f\ : \ \int_{\partial D} \omega f\, d\tau =0\},
$$
where $\Gamma^\nu_\omega$ is considered as an operator on solutions of (\ref{1}).
\end{lemma}

{\sc Remark 1.} The condition $0 \notin \sigma_D(\widetilde{H})$ we understand in the sense that the Dirichlet problem for (\ref{4}) is uniquely solvable. That may be essential if $\omega$ has zeroes, on $D$, implying  singularities of $\widetilde{u}$.

Lemma 2 follows from relation (\ref{14}), the fact that relations (\ref{3}) determine $\theta$ from $\varphi$ up to summands of the form
$\mathrm{const}\cdot \omega^{-1}$ and $\varphi$ from $\theta$ up to summands of the form $\mathrm{const}\cdot \omega$, and from Lemma 1 applied to equation (\ref{4}).

For equation (\ref{1}) and boundary data (\ref{11}) and (\ref{12}) we consider the following
Poincare--Steklov operators
$\Lambda_{u,\omega}$ and $\Lambda^{-1}_{u,\omega}$:
 \begin{equation}
 \label{16}
 \Lambda_{u,\omega}: \Gamma^\tau_\omega \varphi \to \Gamma^\nu_\omega \varphi, \ \ \ 0 \notin \sigma_D(H)
 \end{equation}
 and
 \begin{equation}
 \label{17}
 \Lambda^{-1}_{u,\omega}: \Gamma^\nu_\omega \varphi \to \Gamma^\tau_\omega \varphi, \ \ \ 0 \notin \sigma_D(\widetilde{H}).
 \end{equation}

\begin{proposition}
Assuming that $0 \notin \sigma_D(H)$, the operators $\Phi_u$ and $\Lambda_{u,\omega}$
are related by the following formulas:
$$
\Lambda_{u,\omega}\Gamma^\tau_\omega\varphi = (\Phi_u \omega I \omega^{-1} - \omega_\nu I \omega^{-1})\Gamma^\tau_\omega \varphi,
$$
$$
\Phi_u \varphi = (\Lambda_{u,\omega}\Gamma^\tau_\omega + \omega_\nu\omega^{-1})\varphi,
$$
where $\omega$ (in the cases where $\omega$ is not a lower index), $\omega^{-1}$, and $\omega_\nu$ denote the operators of multiplication by the corresponding functions on $\partial D$ and
$$
I \omega^{-1} f(\tau) = \int^\tau_0 \omega^{-1}(t)f(t)\,dt
$$
for $f \in \mathrm{Im}\, \Gamma^\tau_\omega$.
\end{proposition}

Proposition 1 follows from the definitions given by formulas (\ref{9}), (\ref{11}), and (\ref{16}),
and from straightforward computation.

\begin{proposition}
If $0 \notin \sigma_N(H) \cup \sigma_D(\widetilde{H})$, the operators $\Phi^{-1}_u$ and $\Lambda^{-1}_{u,\omega}$
are related as follows:
$$
\Lambda^{-1}_{u,\omega}\Gamma^\nu_\omega \varphi = \Gamma^\tau_\omega \Phi^{-1}_u
\left(\mathrm{Id} - \frac{\omega_\nu}{\omega}\Phi^{-1}_u \right)^{-1} \Gamma^\nu_\omega \varphi,
$$
where $\mathrm{Id}$ is the identity operator. Therewith,
$$
\mathrm{Ker}\, \left(\mathrm{Id} - \frac{\omega_\nu}{\omega}\Phi^{-1}_u\right) = \{c \omega_\nu \ : \ c\in \C\}
$$
and the inverse operator $\left(\mathrm{Id} - \frac{\omega_\nu}{\omega}\Phi^{-1}_u\right)^{-1}$ is defined on $\mathrm{Im}\, \Gamma^\nu_\omega$
up to summands of the form $c \omega_\nu, c \in \C$.
\end{proposition}

Proposition 2 follows from formulas (\ref{10}), (\ref{11}), (\ref{12}), (\ref{15}), and (\ref{17}), from the representation
$$
\Gamma^\nu_\omega \varphi = \left(\mathrm{Id} - \frac{\omega_\nu}{\omega}\Phi^{-1}_u\right) \varphi_\nu,
$$
and from straightforward computations.

In the assumptions of this section on $D$ and $\omega$ we have the following result.

\begin{theorem}
Let $0 \notin \sigma_D(H) \cup \sigma_D(\widetilde{H})$. Then the following formuals hold:
$$
\Lambda_{\widetilde{u},\omega^{-1}} = - \Lambda^{-1}_{u,\omega},
$$
$$
\Lambda^{-1}_{\widetilde{u},\omega^{-1}} = - \Lambda_{u,\omega}.
$$
\end{theorem}

Theorem 1 follows from relations (\ref{13}) and (\ref{14}) and from Lemmas 1 and 2 which describe the domains of definition of the operators
$\Lambda_{u,\omega}, \Lambda^{-1}_{u,\omega}, \Lambda_{\widetilde{u},\omega^{-1}}$, and
$\Lambda^{-1}_{\widetilde{u},\omega^{-1}}$ under the assumption of Theorem 1.

Propositions 1 and 2 and Theorem 1 provide ways to find the operators $\Phi_{\widetilde{u}}$ and $\Phi^{-1}_{\widetilde{u}}$ from $\Phi_u$ and $\Phi^{-1}_u$ and from the restrictions of $\omega$ and $\omega_\nu$ onto $\partial D$. In this case it is required only the invertibility of operators acting on functions on $\partial D$. In this sense these methods are essentially simpler than the direct reconstruction of $\Phi_{\widetilde{u}}$ and $\Phi^{-1}_{\widetilde{u}}$ from $\widetilde{u}$, when it is necessary to invert an operator acting on functions defined on the whole domain $D$. This effect  becomes
quite obvious for the Darboux transformation that we are demonstrating in the next section.

\section{The action of the Darboux transformation on Poincare-Steklov operators}

Let us consider the Schr\"odinger equation
\begin{equation}
\label{18}
H_1 \psi = \left(-\frac{d^2}{dx^2} + u\right) \psi = E \psi, \ \ \ E=\mu^2
\end{equation}
on the interval $D=]a,b[ \subset \R$
and the Schr\"odinger equation
\begin{equation}
\label{19}
\widetilde{H}_1 \widetilde{\psi} = \left(-\frac{d^2}{dx^2} + \widetilde{u}\right)\widetilde{\psi} = E \widetilde{\psi},
\end{equation}
which is obtained from (\ref{18}) by using the Darboux transformation determined by the solution
$\omega_1$ of the equation
\begin{equation}
\label{20}
H_1 \omega_1 = \kappa^2 \omega_1.
\end{equation}
We assume that $u(x)$ is a regular function of the closed interval $[a,b] = D \cup \partial D$.

For an equation of type (\ref{18}) we consider the operators $Q_u$ and $Q^{-1}_u$
such that
\begin{equation}
\label{21}
Q_u: \psi|_{\partial D} \longrightarrow  \frac{d \psi}{dx}|_{\partial D} ,
\end{equation}
\begin{equation}
\label{22}
Q^{-1}_u:  \frac{d\psi}{dx}|_{\partial D} \longrightarrow  \psi|_{\partial D},
\end{equation}
i.e., the DN and ND operators of the form  (\ref{9}) and (\ref{10}), where, for simplicity, the derivation along the out normal
$\nu$ is replaced by the derivation in $x$.

Let $\sigma_D(H_1)$ and $\sigma_N(H_1)$ define the spectra of the operators defined by $H_1$ 
and the Dirichlet and Neumann
conditions, respectively.

For solutions $\psi$ of (\ref{18}) we also consider the following boundary condition on $\partial D$:
$$
\Gamma^x_{\omega_1}\psi = \omega_1 \left(\frac{\psi}{\omega_1}\right)_x = \psi_x - \frac{\omega_{1,x}}{\omega_1}\psi,
$$
where $\omega_1$ is the given solution of  (\ref{20}) which determines the Darboux trans\-for\-ma\-tion.

We assume that $\omega_1$ has no zeroes on  $\partial D$.

The relations (\ref{6}) and (\ref{7}) on $\partial D$ can be rewritten as follows:
\begin{equation}
\label{23}
\frac{1}{\kappa - \mu}\, \Gamma^x_{\omega_1^{-1}}\widetilde{\psi} = \psi,
\end{equation}
\begin{equation}
\label{24}
\widetilde{\psi} = \frac{1}{\kappa+\mu}\, \Gamma^x_{\omega_1}\psi.
\end{equation}

\begin{lemma}
Assuming that $0 \notin \sigma_D(\widetilde{H}_1)$ and $\kappa \neq \pm \mu$, we have
$$
\mathrm{Ker}\, \Gamma^x_{\omega_1} = 0, \ \ \ \mathrm{Coker}\,\Gamma^x_{\omega_1}=0,
$$
where $\Gamma^x_{\omega_1}$is considered as an operator on solutions of (\ref{18}).
\end{lemma}

{\sc Remark 2.} The condition $0 \notin \sigma_D(\widetilde{H}_1)$ is understood in the sense similar to Remark 1.

Lemma 3 follows from (\ref{7}), (\ref{24}), and the two-dimensionality of the space of functions on $\partial D$ and the spaces of solutions of (\ref{18}) and (\ref{19}).

\begin{theorem}
Assuming that $0 \notin \sigma_N(H_1) \cup \sigma_D(\widetilde{H}_1)$ and $\kappa \neq \pm \mu$,
the following formula is valid:
\begin{equation}
\label{25}
Q_{\widetilde{u}} = -\frac{\omega_{1,x}}{\omega_1} + (\kappa^2-\mu^2)Q^{-1}_u \left(\mathrm{Id} - \frac{\omega_{1,x}}{\omega_1}Q^{-1}_u\right)^{-1},
\end{equation}
where $\frac{\omega_{1,x}}{\omega_1}$ is the operator of multiplication by the corresponding function on $\partial D$.
\end{theorem}

The formula (\ref{25}) explicitly specifies the transformation of the DN and ND operators due to the Darboux transformation.
Such formulas can be used, in particular, for testing the algorithms for reconstructing the potential from measurements at the boundary.

{\sc Proof of Theorem 2.}
We have the representation
\begin{equation}
\label{26}
\Gamma^x_{\omega_1}\psi = \left(\mathrm{Id}-\frac{\omega_{1,x}}{\omega_1}Q^{-1}_u\right)\psi_x.
\end{equation}

It follows from Lemma 3 and (\ref{26}) that the operator
$\mathrm{Id} - \frac{\omega_{1,x}}{\omega_1}Q^{-1}_u$ is invertible.
Further, using (\ref{21}) and (\ref{22}), we rewrite the relations (\ref{23}) and (\ref{24}) in the form
\begin{equation}
\label{27}
\left(Q_{\widetilde{u}} + \frac{\omega_{1,x}}{\omega_1}\right)\widetilde{\psi} = (\kappa-\mu)\psi,
\end{equation}
\begin{equation}
\label{28}
\widetilde{\psi} = \frac{1}{\kappa+\mu}\left(\mathrm{Id} - \frac{\omega_{1,x}}{\omega_1}Q^{-1}_u\right)\psi_x.
\end{equation}
With the help of (\ref{21}) and (\ref{28}), we get
$$
\psi_x = (\kappa+\mu)\left(\mathrm{Id} - \frac{\omega_{1,x}}{\omega_1}Q^{-1}_u\right)^{-1}\widetilde{\psi},
$$
\begin{equation}
\label{29}
\psi= (\kappa+\mu)Q^{-1}_u \left(\mathrm{Id} - \frac{\omega_{1,x}}{\omega_1}Q^{-1}_u\right)^{-1} \widetilde{\psi}.
\end{equation}
Formula (\ref{25}) follows from (\ref{27}) and (\ref{29})
after equating the expressions for $\psi$. Theorem 2 is proved.

{\sc Example.} Let us consider the operator $H = -\frac{d^2}{dx^2}$ on the interval $]a,b[$ with $0<a<b$ and its Darboux transformation determined by $\omega_1=x$. We have
$$
u=0, \ \ \ \widetilde{u} = \frac{2}{x^2}.
$$
Let us take
$$
\mu \neq \kappa,  \ \ \ \mbox{where $\kappa=0$}.
$$
A general solution of the equation
\begin{equation}
\label{30}
H\psi = -\psi^{\prime\prime} = \mu^2 \psi
\end{equation}
has the form
$$
\psi = \alpha \cos \mu x + \beta \sin \mu x.
$$
Therefore,
$$
\left(\begin{array}{c} \psi(a) \\ \psi(b) \end{array}\right) = \left(\begin{array}{cc}\cos \mu a & \sin \mu a \\
\cos \mu b & \sin \mu b \end{array}\right)\left(\begin{array}{c} \alpha \\ \beta \end{array}\right) = A\left(\begin{array}{c} \alpha \\ \beta \end{array}\right),
$$
 $$
\left(\begin{array}{c} \psi^\prime(a) \\ \psi^\prime(b) \end{array}\right) = \mu \left(\begin{array}{cc}-\sin \mu a & \cos \mu a \\
-\sin \mu b & \cos \mu b \end{array}\right)\left(\begin{array}{c} \alpha \\ \beta \end{array}\right) = B\left(\begin{array}{c} \alpha \\ \beta \end{array}\right),
$$
and the DN operator for the problem (\ref{30}) takes the form
$$
Q_0= BA^{-1}  =
\frac{\mu}{\sin \mu (b-a)} \left(\begin{array}{cc} - \cos \mu (b-a) & 1 \\ -1 & \cos \mu (b-a)\end{array}\right).
$$
It is clear that $\sin \mu (b-a) \neq 0$ if and only if $\mu^2$ does not belong to the spectrum of the Dirichlet problem
for $H = -\frac{d^2}{dx^2}$.

The Darboux transformation for solutions takes the form
$$
\widetilde{\psi} = \frac{1}{\mu}\left(\frac{d}{dx} - \frac{1}{x}\right)\psi ,
$$
which implies that
$$
\widetilde{\psi}^{\prime} = \frac{1}{\mu}\left(\psi^{\prime\prime}+\frac{1}{x^2}\psi-\frac{1}{x}\psi^\prime\right) =
\frac{1}{\mu} \left(-\mu^2\psi + \frac{1}{x^2}\psi -\frac{1}{x}\psi^\prime\right)
$$
and the DN operator for the problem
$$
\left(-\frac{d^2}{dx^2} + \frac{2}{x^2}\right)\widetilde{\psi} = \mu^2 \widetilde{\psi}
$$
is equal to
$$
Q_{2/x^2} = \frac{1}{\mu}\left( -\mu^2 + \frac{1}{x}\left(\frac{1}{x}-Q_0\right)\right)\mu\left(Q_0 - \frac{1}{x}\right)^{-1} =
-\frac{1}{x} - \mu^2 \left(Q_0-\frac{1}{x}\right)^{-1},
$$
which is a particular case of (\ref{25}).

\end{document}